Article: CJB/2010/9

# Miquel circles and Cevian lines

**Christopher J Bradley**

**Abstract**: Two theorems are presented concerning the Miquel point configuration, when the operative points on the sides of the triangle are the feet of Cevians,

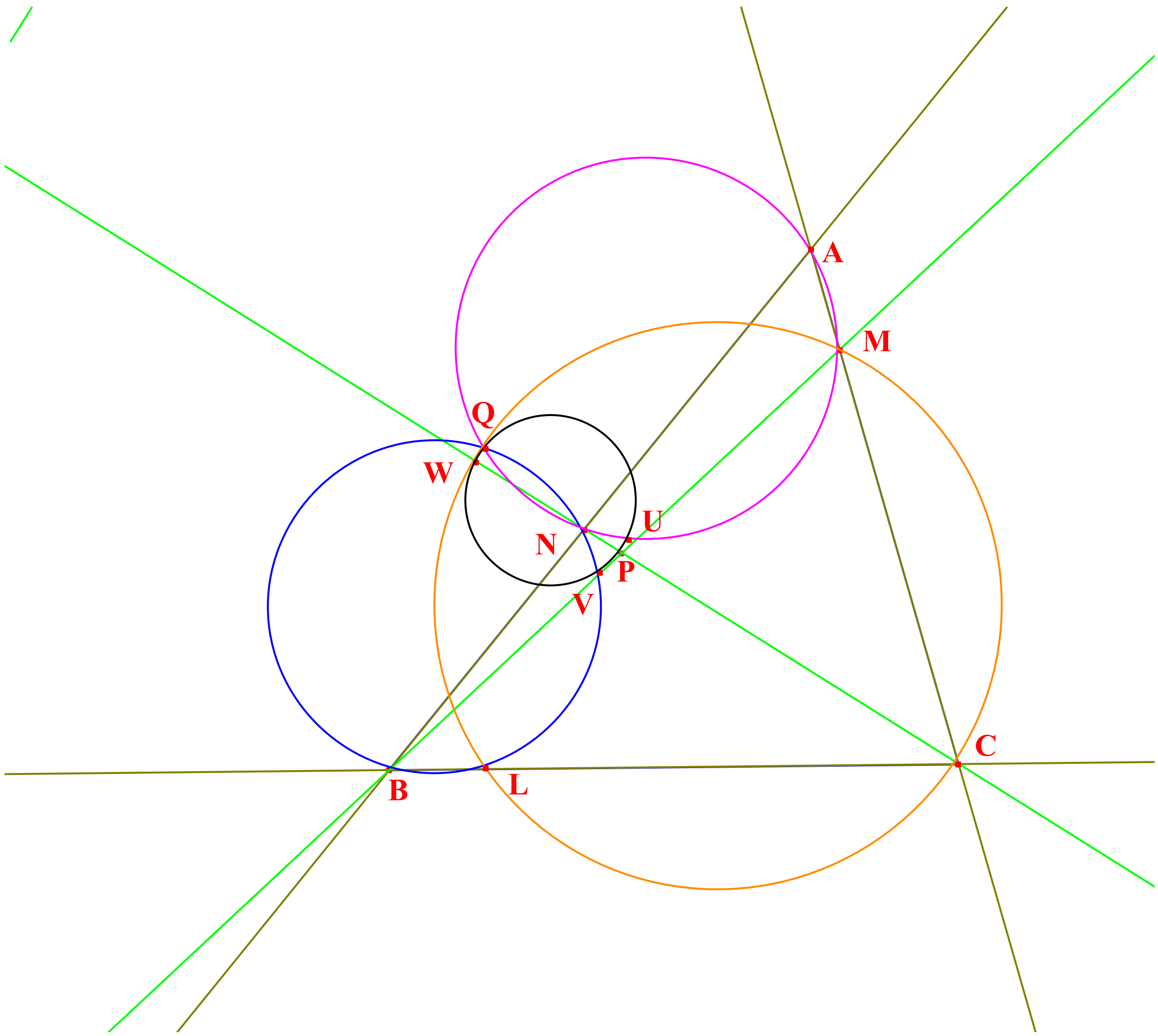

**Fig. 1**

**Miquel point construction from the feet of Cevians**

## 1. Introduction

Given a triangle ABC it is well-known that if you choose L, M, N on the sides BC, CA, AB respectively (but not at the vertices), then circles AMN, BNL, CLM have a common point Q commonly known as the Miquel point [1, 2]. An interesting question is whether there exist any further properties when L, M, N are the feet of Cevian lines through a point P. We answer the question in the affirmative and establish the theorem that both P and Q lie on a circle UVW, where U is the intersection of AP with circle AMN, V is the intersection of BP with circle BNL and W is the intersection of CP with circle CLM. See Fig. 1. We also have a second result which also only holds when L, M, N are the feet of Cevians and this is that circles AQL, BQM, CQN have a common point R. See Figure 2. The proofs given below involve areal co-ordinates. The definition of areal co-ordinates and how to use them may be found in Bradley [3].

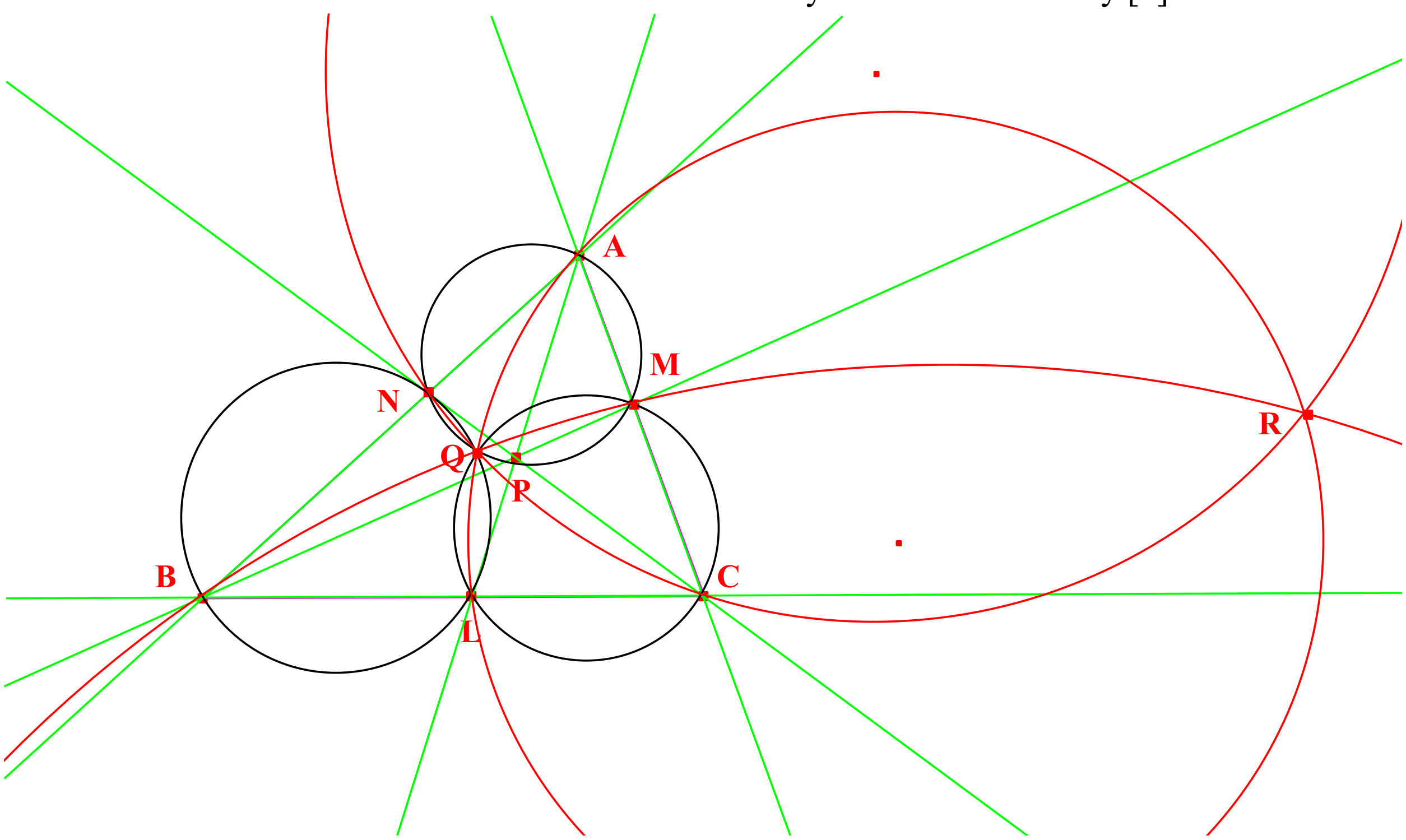

**Fig. 2**

## 2. The Miquel circles

Let P have co-ordinates (l, m, n) so that the co-ordinates of L, M, N are respectively L(0, m, n), M(l, 0, n), N(l, m, 0). The equation of any circle is

$$a^2yz + b^2zx + c^2xy + (x + y + z)(ux + vy + wz) = 0. \qquad (2.1)$$

For circle AMN, we insert the co-ordinates of A, M, N in Equation (2.1) and find three equations for u, v, w, which give $u = 0, v - -\frac{c^2 l}{l+m}, w = -\frac{b^2 l}{n+l}$. We now substitute back in Equation (2.1) to get the equation of AMN, which is

$$c^2\, l(n+\, l)y^2 + b^2 l(l+m)z^2 - \big(a^2(l+m)(n+l) - b^2 l(l+m) - c^2 l(n+l)\big)yz - b^2 n(l+m)zx - c^2 m(n+l)xy = 0. \quad (2.2)$$

The equations of BMN, CNL may now be written down by cyclic change of x, y, z and a, b, c and l, m, n.

## 3. The Miquel Point Q

It may now be checked that these three circles have a common point Q with co-ordinates (x, y, z),
where

$$x = a^2(l+m)(n+l)(-a^2mn(l+m)(n+l) + b^2nl(l+m)(m+n) + c^2lm(n+l)(m+n)),$$

$$y = b^2(l+m)(m+n)(a^2mn(l+m)(n+l) - b^2nl(l+m)(m+n) + c^2lm(n+l)(m+n)),$$

$$z = c^2(m+n)(n+l)\,\big(a^2mn(l+m)(n+l) + b^2nl(l+m)(m+n) - c^2lm(n+l)(m+n)\big). \quad (3.1)$$

## 4. The points U, V, W

The point U is the intersection of AP and circle AMN. AP has equation ny = mz and AMN has Equation (2.2). It follows that U has co-ordinates (x, y, z), where

$$x = -a^2mn(l+m)(n+l) + b^2nl(l+m)(m+n) + c^2lm(m+n)(n+l),$$
$$y = m\big(b^2n^2(l+m) + c^2m^2(n+l)\big), \quad (4.1)$$
$$z = n\big(b^2n^2(l+m) + c^2m^2(n+l)\big).$$

The co-ordinates of V, W now follow by cyclic change of x, y, z and a, b, c and l, m, n.

## 5. The circle UVWPQ

The co-ordinates of U, V, W are now put into an equation of the form (2.1). The result is three equations to determine u, v, w. These are now substituted back and the result is the equation of circle UVW, which is

$$(m+n)\,\big(b^2n^2(l+m) + c^2m^2(n+l)\big)x^2 + \ldots + \cdots \big(l^2(m+n)\big(b^2(l+m) + c^2(n+l)\big) - a^2(l+m)(n+l)(lm+2mn+nl)\big)yz - \cdots - \cdots = 0. \quad (5.1)$$

It may now be verified that both P and Q lie on this circle.

## 6. The circles AQL, BQM, CQN

Starting from Equation (2.1) we may insert the co-ordinates of A, Q, and L to obtain three equations for u, v, w. Substituting these values back into the equation and simplifying we obtain the equation of circle AQL, which is

$$a^2c^2mn(n+l)y^2 - a^2b^2mn(l+m)z^2 + a^2\left(b^2n^2(l+m) - c^2m^2(n+l)\right)yz + b^2\left((m+n)\left(b^2n(l+m) - c^2m(n+l)\right) - a^2mn(l+m)\right)zx + c^2\left((m+n)\left(b^2n(l+m) - c^2m(n+l)\right) + a^2mn(l+m)\right)xy = 0. \tag{6.1}$$

The equations of the circles BQM, CQN may now be written down by cyclic change of x, y, z and a, b, c and l, m, n. Unfortunately the algebra computer package *DERIVE* was not powerful enough to deduce the co-ordinates of the point of intersection of these three circles, other than Q.

## 7. The centres of the three circles

As the three circles are known to have the common point Q, then they will have a second common point R if, and only if, their centres lie on a line. This is because they are then bound to form part of a coaxal system of circles. We refer again to Bradley [3] for how to obtain the centre of a circle. What has to happen is that the polar of the centre is the line at infinity x + y + z = 0. This means that the centre of the conic with equation

$$ux^2 + vy^2 + wz^2 + 2fyz + 2gzx + 2hxy = 0 \tag{7.1}$$

must have co-ordinates X : Y : Z, where

$$X = vw - gv - hw - f^2 + fg + hf, Y = wu - hw - fu - g^2 + gh + fg,$$
$$Z = uv - fu - gv - h^2 + hf + gh. \tag{7.2}$$

Using Equation (6.1) we find the centre of circle AQL has co-ordinates (X, Y, Z), where

$$X = a^2(n+l)\left(a^4n(l+m)(n+l) - a^2\left(b^2n(l+m)(n+2l) + c^2(n+l)(lm+mn+2nl)\right) + l\left(b^4n(l+m) + b^2c^2(lm-2mn-n^2) + c^4(n+l)(m+n)\right)\right)\left(c^2l(m+n) - a^2n(l+m)\right),$$

$$Y = b^2(n+l)\left(a^2n(l+m) - c^2l(m+n)\right)\left(a^4n^2(l+m) - a^2(b^2n^2(l+m) + c^2(l^2(m-n) + ln^2 - mn^2)) + b^2c^2l^2(m+n) - c^4l^2(m+n)\right),$$

(7.3)

$$Z = c^2(n+l)\left(a^2n(l+m) - c^2l(m+n)\right)\left(a^4n(l+m)(n+l) - a^2\left(b^2n(l^2 + 2lm - mn) + c^2(n+l)(lm+mn+2nl)\right) + l(m+n)\left(b^4n - b^2c^2(l+2n) + c^4(n+l)\right)\right).$$

The centres of circles BQM and CQN may now be written down by cyclic change of X, Y, Z and a, b, c and l, m, n. A 3 x 3 determinantal test with the co-ordinates of the centres of these three circles results in a zero value for the determinant, implying that the three centres are collinear, and hence that there must be a second point of intersection R.

Flat 4,
Terrill Court,
12-14 Apsley Road,
BRISTOL BS8 2SP.